# A BERRY–ESSEEN THEOREM FOR FEYNMAN–KAC AND INTERACTING PARTICLE MODELS

BY PIERRE DEL MORAL AND SAMY TINDEL

*Laboratoire Jean Alexandre Dieudonné and Université Henri Poincaré*

In this paper we investigate the speed of convergence of the fluctuations of a general class of Feynman–Kac particle approximation models. We design an original approach based on new Berry–Esseen type estimates for abstract martingale sequences combined with original exponential concentration estimates of interacting processes. These results extend the corresponding statements in the classical theory and apply to a class of branching and genealogical path-particle models arising in nonlinear filtering literature as well as in statistical physics and biology.

**1. Introduction.** Feynman–Kac distribution flows and their particle interpretations arise in the modeling and the numerical solving of a variety of problems including directed polymer simulations in biology and industrial chemistry, nonlinear filtering in advanced signal processing and Bayesian statistics methodology, rare event estimation in telecommunication and computer systems analysis as well as physics in the spectral analysis of Schrödinger operators and in the study of particle absorptions. Their asymptotic behavior as the size of the systems and/or the time parameter tend to infinity has been the subject of various research articles. For more details on both the theoretical and applied aspects of the topic we refer the reader to the review article [5] and references therein.

To better connect this study with existing and related articles in the literature we give a brief discussion on the fluctuation analysis of these models. The first "local" central limit theorems were presented in [1]. These fluctuations were restricted to local sampling errors of an abstract class of genetic type particle model. This study was extended in [2] in the spirit of Shiga and Tanaka's celebrated article [9] to particle and McKean path-measures. This approach to fluctuations in path space was centered around Girsanov type change of measures techniques and a theorem of Dynkin and Mandelbaum on symmetric statistics [7]. This strategy entirely relies on appropriate regularity conditions on the Markov kernels which are not satisfied for genealogical tree evolution models as the ones described in [1]. Another drawback of this approach is that the description of resulting limiting variance is not explicit but expressed in term of the inverse of an $\mathbb{L}_2$ integral operator.

Donsker type theorems and an explicit computation of the limiting variance in terms of Feynman–Kac semigroups were further developed in [4] in the context of particle density profile models. These explicit functional formulations were the starting point of a new approach to central limit theorems based on judicious martingale decompositions and Feynman–Kac semigroup techniques [3, 5].

The main objective of the current article is to complete and further extend these studies with the analysis of the speed of convergence of fluctuations.

The article is organized as follows. In Section 1.1 we describe the Feynman–Kac and the particle models discussed in this article. In Section 1.2 we present our main results and specify the set of regularity conditions needed in the sequel. Section 2 is concerned with a precise Berry–Esseen type estimate for abstract martingale sequences. In Section 3 we show how these martingale fluctuations apply to a sufficiently regular class of McKean particle interpretations.

We end this section with some rather standard and classical notation that will be of current use in the article.

By $\mathcal{M}(E)$ we denote the set of all bounded and positive measures on a measurable space $(E, \mathcal{E})$, by $\mathcal{P}(E) \subset \mathcal{M}(E)$ we denote the subset of probability measures on $(E, \mathcal{E})$ and by $\mathcal{B}_b(E)$ the Banach space of all bounded $\mathcal{E}$-measurable functions $f$ on $E$ equipped with the uniform norm $\|f\| = \sup_{x \in E} |f(x)|$. We also let $\mathrm{Osc}(E) \subset \mathcal{B}_b(E)$ be the subset of all bounded measurable functions with oscillations $\mathrm{osc}(f) = \sup_{(x,y)} |f(x) - f(y)| \leq 1$.

We finally recall that a bounded and positive integral operator $Q$ from $(E, \mathcal{E})$ into another measurable space $(E', \mathcal{E}')$ generates two operators, one acting on functions $f \in \mathcal{B}_b(E')$ and taking values in $\mathcal{B}_b(E)$, the other acting on measures $\mu \in \mathcal{M}(E)$ into $\mathcal{M}(E')$ and defined by

$$Q(f)(x) = \int_{E'} Q(x, dx') f(x'), \qquad \mu Q(dx') = \int_E \mu(dx) Q(x, dx').$$

To clarify the presentation we shall slightly abuse the notation, and we often write $Q(f - Q(f))^2$ for the function $x \to Q(f - Q(f)(x))^2(x)$.

Finally, we shall use the letter $c$ to denote any nonnegative and universal constant whose values may vary from line to line but which does not depend on the time parameter nor on the Feynman–Kac models.

1.1. *Description of the models.* We consider some collections of measurable spaces $(E_n, \mathcal{E}_n)_{n \in \mathbb{N}}$, of Markov transitions $M_{n+1}(x_n, dx_{n+1})$ from $E_n$ to $E_{n+1}$, and bounded $\mathcal{E}_n$-measurable and strictly positive functions $G_n$ on $E_n$. We assume that the latter are chosen so that for any $n \in \mathbb{N}$ we have

$$r_n = \sup_{(x_n, y_n) \in E_n^2} (G_n(x_n)/G_n(y_n)) < \infty. \tag{1}$$

We associate to the pair $(G_n, M_n)$ the Boltzmann–Gibbs transformation $\Psi_n$ on $\mathcal{P}(E_n)$ and the mapping $\Phi_{n+1}$ from $\mathcal{P}(E_n)$ into $\mathcal{P}(E_{n+1})$ given for any $(x_n, \mu_n) \in (E_n, \mathcal{P}(E_n))$ by

$$\Psi_n(\mu_n)(dx_n) = G_n(x_n)\mu_n(dx_n)/\mu_n(G_n),$$

$$\Phi_{n+1}(\mu_n) = \Psi_n(\mu_n)M_{n+1}.$$

For any $\eta_0 \in \mathcal{P}(E_0)$ we denote by $\mathbb{E}_{\eta_0}(\cdot)$ the expectation operator with respect to the distribution of a Markov chain $X_n$ with initial distribution $\eta_0$ and elementary transitions $M_n$. We consider the distribution flow $\eta_n \in \mathcal{P}(E_n)$, $n \in \mathbb{N}$, also called *Feynman–Kac flow* in the sequel, defined for any $f_n \in \mathcal{B}_b(E_n)$ by the Feynman–Kac formulae

$$\eta_n(f_n) = \gamma_n(f_n)/\gamma_n(1) \qquad \text{with } \gamma_n(f_n) = \mathbb{E}_{\eta_0}\left[f_n(X_n) \prod_{0 \leq p < n} G_p(X_p)\right], \tag{2}$$

with the convention $\prod_{\varnothing} = 1$. Using the multiplicative structure of the Feynman–Kac model and the Markov property, one readily checks that the flow $\eta_n$ satisfies the nonlinear equation

$$\eta_{n+1} = \eta_n K_{n+1, \eta_n}, \tag{3}$$

where $(K_{n+1,\mu_n})_{n \in \mathbb{N}, \mu_n \in \mathcal{P}(E_n)}$ is a *nonunique* collection of Markov transitions satisfying the compatibility condition

$$\forall n \in \mathbb{N} \ \forall \mu_n \in \mathcal{P}(E_n) \qquad \mu_n K_{n+1, \mu_n} = \Phi_{n+1}(\mu_n). \tag{4}$$

These collections of transitions are often called the McKean interpretations of (3). Notice that the compatibility relation (4) is satisfied if we take

$$K_{n+1,\mu_n}(x_n, \cdot) = \varepsilon_n G_n(x_n) M_{n+1}(x_n, \cdot) + (1 - \varepsilon_n G_n(x_n))\Phi_{n+1}(\mu_n) \tag{5}$$

for any nonnegative constant $\varepsilon_n$ such that $\varepsilon_n G_n(x_n) \in [0, 1]$. We finally notice that the random variables $X_n$ may represent the path of an auxiliary Markov chain $X'_p$ from the origin up to time $n$ and taking values in some Hausdorff topological spaces $E'_p$; that is, we have

$$X_n = (X'_0, \ldots, X'_n) \in E_n = (E'_0 \times \cdots \times E'_n). \tag{6}$$

For each $N \geq 1$ we denote by $m_N$ the mapping from the product space $E^N$ into $\mathcal{P}(E)$ which associates to each configuration $x = (x^i)_{1 \leq i \leq N} \in E^N$ the empirical measure $m_N(x) = \frac{1}{N} \sum_{i=1} \delta_{x^i}$. The interacting particle system associated to a given McKean interpretation is defined as a Markov chain $\xi_n^{(N)} = (\xi_n^{(N,i)})_{1 \leq i \leq N}$ taking values in the product spaces $E_n^N$ with initial distribution $\eta_0^{\otimes N}$ and elementary transitions

$$\text{Prob}(\xi_n^{(N)} \in dx_n | \xi_{n-1}^{(N)}) = \prod_{i=1}^{N} K_{n, m_N(\xi_{n-1}^{(N)})}(\xi_{n-1}^{(N,i)}, dx_n^i), \tag{7}$$

where $dx_n = \times_{1 \leq i \leq N} dx_n^i$ stands for an infinitesimal neighborhood of the point $x_n = (x_n^i)_{1 \leq i \leq N} \in E_n^N$.

Under appropriate regularity conditions on the McKean transitions kernels $K_{n,\mu_n}$ it is known that in some sense the particle measures

$$\eta_n^N = m_N(\xi_n^{(N)})$$

converge as $N$ tends to infinity to the desired distributions $\eta_n$.

To illustrate this model we note that the particle interpretation of the Feynman–Kac flow associated to McKean transitions (5) forms a two-step selection/mutation genetic algorithm. The particular situation where $\varepsilon_n = 0$ corresponds to a simple genetic model with an overlapping mutation/selection transition. In the same vein the corresponding particle interpretation model of the Feynman–Kac path measures associated to the chain (6) forms is a genetic type algorithm taking values in path space. Note that in this situation the path-particles have the form

$$\xi_n^{(N,i)} = (\xi_{0,n}^{(N,i)}, \xi_{1,n}^{(N,i)}, \ldots, \xi_{n,n}^{(N,i)}) \in E_n = (E_0' \times \cdots \times E_n').$$

In addition, if the potential functions only depend on terminal values in the sense that $G_n(x_0', \ldots, x_n') = G_n'(x_n')$ for some potential function $G_n'$ on $E_n'$, then the resulting path-particle model can be interpreted as a genealogical tree evolution model.

As traditionally, to clarify the presentation we slightly abuse the notation, by suppressing the size index $N$, and we write $(m(x), \xi_n, \xi_n^i)$ instead of $(m_N(x), \xi_n^{(N)}, \xi_n^{(N,i)})$.

1.2. *Statement of some results.* For any sequence of $\mathcal{F}_n^N$-adapted random variables $Z_n^N$ defined on some filtered probability spaces $(\Omega^N, (\mathcal{F}_n^N)_{n \geq 0}, \mathbb{P}^N)$, we denote by $\Delta Z_n^N$ the difference process $\Delta Z_n^N = Z_n^N - Z_{n-1}^N$, with the convention $\Delta Z_0^N = Z_0^N$ for $n = 0$. If $\Delta M_n^N$ is a given $\mathcal{F}_n^N$-martingale difference, then we denote by $M_n^N$ the $\mathcal{F}_n^N$-martingale defined by $M_n^N = \sum_{p=0}^{n} \Delta M_p^N$. We recall that its increasing process $\langle M^N \rangle_n$ is given by

$$\langle M^N \rangle_n = \sum_{p=0}^{n} \mathbb{E}^N[(\Delta M_p^N)^2 | \mathcal{F}_{p-1}^N],$$

with the convention $\mathcal{F}_{-1}^N = \{\varnothing, \Omega^N\}$ for $p = 0$. It is also convenient to introduce the increasing process $C_n^N = N \langle M^N \rangle_n$ of the normalized martingale $L_n^N = \sqrt{N} M_n^N$.

The example we have in mind is the situation where $(\mathcal{F}_n^N)_{n \geq 0}$ is the natural filtration associated to the particle model (7) and the $\mathcal{F}^N$-martingale difference $\Delta M_n^N = \Delta \mathcal{M}_n^N(f_n)$, with $f_n \in \mathcal{B}_b(E_n)$, is given by the particle $n$th sampling error

$$(8) \qquad \Delta \mathcal{M}_n^N(f_n) = \eta_n^N(f_n) - \eta_{n-1}^N K_{n,\eta_{n-1}^N}(f_n)$$

with the convention $\eta_{-1}^N K_{0,\eta_{-1}^N} = \eta_0$ for $n = 0$. In this situation the increasing processes $\langle \mathcal{M}^N(f) \rangle_n$ and $\mathcal{C}_n^N(f) = \langle \mathcal{L}^N(f) \rangle_n$ of the corresponding martingales $\mathcal{M}_n^N(f)$ and $\mathcal{L}_n^N(f) = \sqrt{N} \mathcal{M}_n^N(f)$ are connected by the formula

$$(9) \qquad \mathcal{C}_n^N(f) = N \langle \mathcal{M}^N(f) \rangle_n = \sum_{p=0}^n \eta_{p-1}^N K_{p,\eta_{p-1}^N}(f_p - K_{p,\eta_{p-1}^N}(f_p))^2.$$

Our first main result concerns a Berry–Esseen theorem for an abstract class of martingale sequences under the following set of conditions:

(H1) For any $n \geq 0$ there exist some constants $a_1(n) < \infty$ and $0 < c_1(n) \leq 1$ such that for any $n \geq 0$ and $\lambda^3 \leq c_1(n) N^{1/2}$ we have, $\mathbb{P}^N$ almost surely,

$$|\mathbb{E}[e^{i\lambda N^{1/2} \Delta M_n^N + (\lambda^2/2) \Delta C_n^N} | \mathcal{F}_{n-1}^N] - 1| \leq a_1(n) \lambda^3 / N^{1/2}.$$

(H2) For any $n \geq 0$ there exists some finite constant $a_2(n) < \infty$ such that for any $N \geq 1$, $\lambda > 0$ and $n \geq 0$,

$$|\mathbb{E}[e^{i\lambda N^{1/2} M_n^N}]| \leq \mathbb{E}[e^{-(\lambda^2/2)\Delta C_n^N}] e^{\lambda^3 a_2(n)/N^{1/2}}.$$

(H3) There exists a nonnegative and strictly increasing deterministic process $C = (C_n)_{n \geq 0}$ as well as some finite constants $0 < a_3(n) < \infty$ such that for any $\varepsilon > 0$ we have

$$\mathbb{E}[e^{\varepsilon N^{1/2} |\Delta C_n^N - \Delta C_n|}] \leq (1 + \varepsilon a_3(n)) e^{\varepsilon^2 a_3^2(n)}.$$

THEOREM 1.1. *Let $M^N = (M_n^N)_{n \geq 0}$ be a sequence of $\mathcal{F}^N$-martingales satisfying conditions* (H1)–(H3) *for some nonnegative and strictly increasing process $C_n$. We let $F_n^N$, respectively, $F_n$, be the distribution function of the random variable $N^{1/2} M_n^N$, respectively, the one of a centered Gaussian random variable with variance $C_n$. Then for any $n \geq 0$ we have*

$$\limsup_{N \to \infty} N^{1/2} \| F_n^N - F_n \| < \infty.$$

Theorem 1.1 does not apply directly to the particle martingale sequence introduced in (8). The first two conditions (H1) and (H2) are rather standard. They can be checked for any kind of any McKean interpretation model using simple and rather standard asymptotic expansions of characteristic functions. The third condition is an exponential continuity condition of the increasing processes introduced in (9). Next we provide a sufficient regularity condition which can be easily checked in various McKean interpretation models. If we set for any $\mu_n \in \mathcal{P}(E_n)$

$$\mu_n^{-1}(0) = \{h \in \mathrm{Osc}(E_n) : \mu_n(h) = 0\},$$

then this condition reads

(H) There exists a collection of uniformly bounded positive measures $\Gamma_{n,f}$ and $\Gamma'_{n+1,f}$ on the sets $\eta_n^{-1}(0)$ and $\eta_{n+1}^{-1}(0)$ and indexed by $n \in \mathbb{N}$ and $f \in \mathrm{Osc}(E_{n+1})$ and such that

$$\|K_{n+1,\mu_n}(f) - K_{n+1,\eta_n}(f)\|$$
$$\leq \int |\mu_n(h)| \Gamma_{n,f}(dh) + \int |\Phi_{n+1}(\mu_n)(h)| \Gamma'_{n,f}(dh).$$

When condition (H) is met we denote by $\Gamma$ the supremum of the total mass quantities $\Gamma'_{n,f}(1)$ and $\Gamma_{n,f}(1)$.

Note that (H) is related to some Lipschitz type regularity of the increasing process and it is clearly met for the McKean transitions given in (5), since we have in this case

$$K_{n+1,\mu_n}(f) - K_{n+1,\eta_n}(f) = (1 - \varepsilon_n G_n)[\Phi_{n+1}(\mu_n) - \Phi_{n+1}(\eta_n)](f).$$

Thus, in this situation, we have that (H) is met with $\Gamma_{n,f} = 0$ and $\Gamma'_{n,f} = \delta_h$ with $h = [f - \eta_{n+1}(f)]$ so that $\Gamma_{n,f}(1) = 0$ and $\Gamma'_{n,f}(1) = 1$. When the parameter $\varepsilon_n = \varepsilon_n(\mu_n)$ in (5) depends on the index measure $\mu_n$, we also find that (H) is met as soon as we have

$$|\varepsilon_n(\mu_n) - \varepsilon_n(\eta_n)| \leq \int |\mu_n(h)| \Lambda_n(dh)$$

for some uniformly bounded positive measures $\Lambda_n$ on $\eta_n^{-1}(0)$.

REMARK 1.2. The above considerations show that condition (H3) is in fact easily verified in most of the classical applications of McKean models, and in particular in the case of nonlinear filtering, for which we refer to [5], for sake of conciseness.

To check the exponential estimates stated in condition (H3) we shall use a refined version of Burkholder's inequality recently presented by the first author with Miclo and Viens in [6]. Roughly speaking, these sharp $\mathbb{L}_p$-estimates combined with some judicious error decompositions lead to the desired exponential concentrations estimates for the normalized sampling error martingale $\mathcal{L}_n^N(f)$ defined by (9) with the limiting increasing process

$$(10) \qquad \mathcal{C}_n(f) = \sum_{p=0}^{n} \eta_{p-1} K_{p,\eta_{p-1}}(f_p - K_{p,\eta_{p-1}}(f_p))^2.$$

Observe that, even if this strategy led to the desired Berry–Esseen estimates on $\mathcal{M}_n^N(f)$, we would still need to transfer these rates of fluctuations to the random field sequence defined by $\sqrt{N}(\eta_n^N - \eta_n)$. One of the most elegant approaches is probably to follow the semigroup techniques and the martingale decompositions developed in [5]. To describe these decompositions with some precision we let $Q_{p,n}$ be the linear Feynman–Kac semigroup associated to the flow $\gamma_n$. To be more specific, we define the semigroup $Q_{p,n}$ by the relation

$$\gamma_n = \gamma_p Q_{p,n},$$

and we associate to $Q_{p,n}$ a "normalized" semigroup $\overline{Q}_{p,n}$, defined for $f_n \in \mathcal{B}_b(E_n)$ by

$$(11) \qquad \overline{Q}_{p,n}(f_n) = \frac{Q_{p,n}(f_n)}{\eta_p Q_{p,n}(1)} = \frac{\gamma_p(1)}{\gamma_n(1)} Q_{p,n}(f_n).$$

If we let $(W_{p,n}^N(f_n))_{p \leq n}$, $f_n \in \mathcal{B}_b(E_n)$, be the random field sequence defined by

$$(12) \quad W_{p,n}^N(f_n) = \sqrt{N}(\eta_p^N - \eta_p)(f_{p,n}) \qquad \text{with } f_{p,n} = \overline{Q}_{p,n}(f_n - \eta_n f_n),$$

then we have the Doob type decomposition

$$W_{p,n}^N(f_n) = \mathcal{B}_{p,n}^N(f_n) + \mathcal{L}_{p,n}^N(f_n),$$

with the predictable and martingale sequences given by

$$\Delta \mathcal{B}_{p,n}^N(f_n) = \sqrt{N}[1 - \eta_{p-1}^N(G_{p-1})/\eta_{p-1}(G_{p-1})]$$
$$\times [\Phi_p(\eta_{p-1}^N)(f_{p,n}) - \Phi_p(\eta_{p-1})(f_{p,n})],$$
$$\Delta \mathcal{L}_{p,n}^N(f_n) = \sqrt{N}[\eta_p^N(f_{p,n}) - \eta_{p-1}^N K_{p,\eta_{p-1}^N} f_{p,n}].$$

The above decomposition is now more or less standard. For the convenience of the reader its proof is housed in the Appendix.

Intuitively speaking, we see from the quadratic structure of the predictable term that it should not influence the fluctuation rate. We will make precise this observation with a Stein type approximation lemma and we will prove the following theorem:

THEOREM 1.3. *Let $f_n \in \mathcal{B}_b(E_n)$ and let $W_{p,n}^N$ be the quantity defined by* (12). *For any McKean interpretation model satisfying condition* (H), *we have*

$$\limsup_{N \to \infty} \sqrt{N} \sup_{u \in \mathbb{R}} \left| \mathbb{P}(W_{n,n}^N(f_n) \leq u\sigma_n(f)) - \int_{-\infty}^{u} e^{-v^2/2} \frac{dv}{(2\pi)^{1/2}} \right| < \infty,$$

*for any $f_n \in \mathcal{B}_b(E_n)$ and $n \geq 0$ with*

$$\sigma_n^2(f) = \sum_{p=0}^{n} \eta_{p-1} K_{p,\eta_{p-1}}(f_{p,n} - K_{p,\eta_{p-1}}(f_{p,n}))^2.$$

**2. An estimate for martingale sequences.** The central limit theorem for sequences of random variables is usually obtained by convergence of characteristic functions. Unsurprisingly, the natural question of determining the speed of convergence in the CLT can also be handled through characteristic functions considerations. The formalization of this idea is due to Berry and Esseen, and can be summarized in the following theorem:

THEOREM 2.1 (Berry–Esseen). *Let $(F_1, F_2)$ be a pair of distribution functions with characteristic functions $(f_1, f_2)$. Also assume that $F_2$ has a derivative with $\|\frac{\partial F_2}{\partial x}\| < \infty$. Then for any $a > 0$ we have*

$$\|F_1 - F_2\| \leq \frac{2}{\pi} \int_0^a \frac{|f_1(x) - f_2(x)|}{x} dx + \frac{24}{a\pi} \left\| \frac{\partial F_2}{\partial x} \right\|.$$

In this section, we will try to apply this theorem to a sequence of martingales satisfying the general set of hypotheses (H1)–(H3) in order to get a sharp asymptotic result for its convergence towards a Gaussian martingale. In order to prepare for the proof of Theorem 1.1 we start with the following key technical lemma.

LEMMA 2.2. *Suppose we are given a sequence of martingales $M^N = (M_n^N)_{n \geq 0}$ with respect to some filtrations $\mathcal{F}_n^N$, satisfying the conditions* (H1)–(H3). *Then, for any $n \geq 0$, there exist a finite constant $a(n) < \infty$, a positive constant $b(n)$ and some $N(n) \geq 1$ such that for any $N \geq N(n)$ and $0 < \lambda \leq b(n)N^{1/2}$,*

$$|\mathbb{E}[e^{i\lambda N^{1/2} M_n^N}] - e^{-(\lambda^2/2)C_n}| \leq a(n) e^{-(\lambda^2/4)\Delta C_n} \frac{\lambda^2(1+\lambda)}{N^{1/2}}.$$

Since the proof of Theorem 1.1 is a simple consequence of the above lemma we have chosen to give it first.

PROOF OF THEOREM 1.1. By Theorem 2.1 and Lemma 2.2 we have, for any $N \geq N(n)$,

$$N^{1/2}\|F_n^N - F_n\|$$
$$\leq \frac{2a(n)}{\pi} \int_0^{b(n)N^{1/2}} e^{-(\lambda^2/4)\Delta C_n} \lambda(1+\lambda)\,d\lambda + \frac{24}{b(n)(2e\pi^3 C_n)^{1/2}}$$
$$\leq \frac{2a(n)}{\pi} \int_0^\infty e^{-(\lambda^2/4)\Delta C_n} \lambda(1+\lambda)\,d\lambda + \frac{24}{b(n)C_n^{1/2}},$$

for some $N(n) \geq 1$ and some finite positive constant $0 < b(n) < \infty$. Invoking the fact that $\Delta C_n > 0$, this ends the proof of the theorem. $\square$

We now come to the proof of Lemma 2.2.

PROOF OF LEMMA 2.2. Let $I_n^N$ be the function defined for any $\lambda \geq 0$ by

$$I_n^N(\lambda) = \mathbb{E}[e^{i\lambda N^{1/2} M_n^N + (\lambda^2/2)C_n}] - 1,$$

and notice that

(13) $$\mathbb{E}[e^{i\lambda N^{1/2} M_n^N}] - e^{-(\lambda^2/2)C_n} = e^{-(\lambda^2/2)C_n} I_n^N(\lambda).$$

Furthermore, we have the easily verified recursive equations

$$I_n^N(\lambda) - I_{n-1}^N(\lambda)$$
$$= \mathbb{E}[e^{i\lambda N^{1/2} M_{n-1}^N + (\lambda^2/2)C_{n-1}}(\mathbb{E}[e^{i\lambda N^{1/2}\Delta M_n^N + (\lambda^2/2)\Delta C_n}|F_{n-1}^N] - 1)],$$

and hence

$$I_n^N(\lambda) - I_{n-1}^N(\lambda) = A(\lambda) + B(\lambda),$$

with

$$A(\lambda) = \mathbb{E}[e^{i\lambda N^{1/2} M_{n-1}^N + (\lambda^2/2)C_{n-1}}(e^{(\lambda^2/2)(\Delta C_n - \Delta C_n^N)} - 1)]$$

and

$$B(\lambda) = \mathbb{E}[e^{i\lambda N^{1/2} M_{n-1}^N + (\lambda^2/2)C_{n-1}}(e^{(\lambda^2/2)(\Delta C_n - \Delta C_n^N)}) \\ \times (\mathbb{E}[e^{i\lambda N^{1/2}\Delta M_n^N + (\lambda^2/2)\Delta C_n^N}|\mathcal{F}_{n-1}^N] - 1)].$$

Using this we obtain

$$|I_n^N(\lambda) - I_{n-1}^N(\lambda)| \leq e^{(\lambda^2/2)C_{n-1}}(A_1(\lambda) + B_1(\lambda)),$$

where

$$A_1(\lambda) = \mathbb{E}[e^{(\lambda^2/2)|\Delta C_n - \Delta C_n^N|} - 1]$$

and
$$B_1(\lambda) = \mathbb{E}[|\mathbb{E}[e^{i\lambda N^{1/2}\Delta M_n^N + (\lambda^2/2)\Delta C_n^N}|\mathcal{F}_{n-1}^N] - 1|e^{(\lambda^2/2)(\Delta C_n - \Delta C_n^N)}].$$

Now, under conditions (H1) and (H3) applied for $\varepsilon = \frac{\lambda^2}{2N^{1/2}}$, we find that

$$|I_n^N(\lambda) - I_{n-1}^N(\lambda)|$$
$$\leq e^{(\lambda^2/2)C_{n-1}}\left[\frac{a_1(n)\lambda^3}{N^{1/2}}\left(1 + \frac{\lambda^2 a_3(n)}{2N^{1/2}}\right)e^{(\lambda^4/4N)a_3^2(n)}\right.$$
$$\left.+ \left(1 + \frac{\lambda^2 a_3(n)}{2N^{1/2}}\right)e^{(\lambda^4/4N)a_3^2(n)} - 1\right]$$
$$= e^{(\lambda^2/2)C_{n-1}}\left[\frac{a_1(n)\lambda^3}{N^{1/2}}\left(1 + \frac{\lambda^2 a_3(n)}{2N^{1/2}}\right)e^{(\lambda^4/4N)a_3^2(n)}\right.$$
$$\left.+ (e^{(\lambda^4/4N)a_3^2(n)} - 1) + \frac{\lambda^2 a_3(n)}{2N^{1/2}}e^{(\lambda^4/4N)a_3^2(n)}\right],$$

for any $0 < \lambda^3 \leq c_1(n)N^{1/2}$. Since for these pairs of parameters $(\lambda, N)$ we have $\lambda^2 \leq N^{1/2}$ (and therefore $\lambda^4 \leq N$), we find that

$$N^{1/2}|I_n^N(\lambda) - I_{n-1}^N(\lambda)| \leq d(n)e^{(\lambda^2/2)C_{n-1}}\lambda^2(1+\lambda),$$

for some finite constant $d(n)$, whose value only depends on $a_i(n)$, $i = 1, 3$, and such that

$$d(n) \leq c e^{a_3^2(n)/4}(1 \vee a_1(n) \vee a_3(n))^2.$$

If we set

$$c_\star(n) = \bigwedge_{p=0}^{n} c_1(p) \ (\leq 1) \quad \text{and} \quad d^\star(n) = \bigvee_{p=0}^{n} d(p),$$

then for any $0 \leq p \leq n$ and any $0 < \lambda^3 \leq c_\star(n)N^{1/2}$, we have that

$$N^{1/2}|I_p^N(\lambda) - I_{p-1}^N(\lambda)| \leq d^\star(n)e^{(\lambda^2/2)C_{n-1}}\lambda^2(1+\lambda).$$

It is now easily verified from these estimates that

$$N^{1/2}|I_n^N(\lambda)| \leq (n+1)d^\star(n)e^{(\lambda^2/2)C_{n-1}}\lambda^2(1+\lambda),$$

from which we conclude that, for any $0 < \lambda^3 \leq c_\star(n)N^{1/2}$,

$$(14) \quad |\mathbb{E}[e^{i\lambda N^{1/2}M_n^N}] - e^{-(\lambda^2/2)C_n}| \leq (n+1)d^\star(n)e^{-(\lambda^2/2)\Delta C_n}\frac{\lambda^2(1+\lambda)}{N^{1/2}}.$$

On the other hand, we have, for any pair $(\lambda, N)$,

$$(15) \quad |\mathbb{E}[e^{i\lambda N^{1/2}M_n^N}] - e^{-(\lambda^2/2)C_n}| \leq |\mathbb{E}[e^{i\lambda N^{1/2}M_n^N}]| + e^{-(\lambda^2/2)C_n},$$

and under condition (H2),

$$|\mathbb{E}[e^{i\lambda N^{1/2}M_n^N}]| \leq \mathbb{E}[e^{-(\lambda^2/2)\Delta C_n^N}]e^{\lambda^3 a_2(n)N^{-1/2}}.$$

Using again (H3) we also find that

$$|\mathbb{E}[e^{i\lambda N^{1/2}M_n^N}]|$$
$$\leq e^{-(\lambda^2/2)\Delta C_n}\left(1 + \frac{\lambda^2 a_3(n)}{2N^{1/2}}\right)e^{(\lambda^4 a_3^2(n))/(4N)}e^{(\lambda^3 a_2(n))/N^{1/2}}$$
$$= \left(1 + \frac{\lambda^2 a_3(n)}{2N^{1/2}}\right)\exp\left[-\frac{\lambda^2}{2}\left(\Delta C_n - \frac{\lambda}{N^{1/2}}\left(2a_2(n) + a_3^2(n)\frac{\lambda}{2N^{1/2}}\right)\right)\right].$$

Recall that $\Delta C_n > 0$, and observe that for any pair $(\lambda, N)$ such that

$$\lambda \leq c^\star(n)N^{1/2} \qquad \text{with } c^\star(n) = [2a_3^{-2}(n) \wedge (2^{-1}\Delta C_n(1 + 2a_2(n))^{-1})],$$

we have

$$\frac{\lambda}{N^{1/2}}\left(2a_2(n) + a_3^2(n)\frac{\lambda}{2N^{1/2}}\right) \leq \frac{\lambda}{N^{1/2}}(2a_2(n) + 1)$$
$$\leq \frac{\Delta C_n}{2}.$$

This yields

$$(16) \qquad |\mathbb{E}[e^{i\lambda N^{1/2}M_n^N}]| \leq \left(1 \vee \frac{a_3(n)}{2}\right)\left(1 + \frac{\lambda^2}{N^{1/2}}\right)e^{-\lambda^2 \Delta C_n/4},$$

and hence, by (15), and for any $\lambda \leq c^\star(n)N^{1/2}$, we find

$$(17) \quad |\mathbb{E}[e^{i\lambda N^{1/2}M_n^N}] - e^{-(\lambda^2/2)C_n}| \leq e^{-(\lambda^2/4)\Delta C_n}(2 \vee a_3(n))\left(1 + \frac{\lambda^2}{N^{1/2}}\right).$$

To take the final step we observe that for any

$$N \geq c_\star(n)/c^\star(n)^3 \quad \text{and} \quad c_\star^{1/3}(n)N^{1/6} \leq \lambda \leq c^\star(n)N^{1/2},$$

we have $1 = c_\star(n)/c_\star(n) \leq c_\star^{-1}(n)\lambda^3/N^{1/2}$, and by (17),

$$(18) \quad |\mathbb{E}[e^{i\lambda N^{1/2}M_n^N}] - e^{-(\lambda^2/2)C_n}| \leq c_\star^{-1}(n)(2 \vee a_3(n))\frac{\lambda^2(1 + \lambda)}{N^{1/2}}e^{-\lambda^2 \Delta C_n/4}.$$

In conjunction with (14) we conclude that for any $N \geq N(n) = c_\star(n)/c^\star(n)^3$ and any $\lambda \leq c^\star(n)N^{1/2}$,

$$|\mathbb{E}[e^{i\lambda N^{1/2}M_n^N}] - e^{-(\lambda^2/4)C_n}| \leq a(n)\frac{\lambda^2}{N^{1/2}}(1 + \lambda)e^{-(\lambda^2/4)\Delta C_n},$$

with $a(n) = [(n + 1)d^\star(n)] \vee [c_\star^{-1}(n)(2 \vee a_3(n))]$. This ends the proof of the lemma. $\square$

**3. Application to interacting processes.** In this section, we prove that Theorem 1.1 can be applied to our particle approximations. We shall go through a series of preliminary results leading to the proof of Theorem 1.3.

The first step is of course to provide some exponential estimates for the particle density profiles. In the next pivotal lemma we describe an original exponential concentration result in terms of the following pair of parameters:

$$\beta(P_{p,n}) = \sup_{(x_p,y_p)\in E_p^2} \|P_{p,n}(x_p,\cdot) - P_{p,n}(y_p,\cdot)\|_{\mathrm{tv}},$$

(19)

$$r_{p,n} = \sup_{(x_p,y_p)\in E_p^2} Q_{p,n}(1)(x_p)/Q_{p,n}(1)(y_p),$$

where $\|\mu - \nu\|_{\mathrm{tv}} = \sup_{A\in\mathcal{E}} |\mu(A) - \nu(A)|$ represents the total variation distance between probabilities and $P_{p,n}$ denotes the Markov transition from $E_p$ into $E_n$ defined by

$$P_{p,n}(x_p, dx_n) = Q_{p,n}(x_p, dx_n)/Q_{p,n}(x_p, E_n).$$

LEMMA 3.1. *For any McKean model we have for every $n \geq 0$, $f_n \in \mathrm{Osc}_1(E_n)$ and $\varepsilon > 0$,*

$$\mathbb{E}[e^{\varepsilon N^{1/2}|\eta_n^N(f_n) - \eta_n(f_n)|}] \leq (1 + \varepsilon 2^{-1/2} b(n)) e^{(\varepsilon b(n))^2/2}$$

*for some finite constant $b(n)$ such that $b(n) \leq 2\sum_{q=0}^{n} r_{q,n}\beta(P_{q,n})$.*

REMARK 3.2. The quantities $(r_{p,n}, \beta(P_{p,n}))$ play an important role in the asymptotic and long time behavior of Feynman–Kac particle approximation models. The above lemma combined with the semigroups approach developed in [5] readily yields uniform exponential concentration properties. To be more specific, let us suppose that $r = \bigvee_n r_n < \infty$. Also assume that there exist some integer parameter $m \geq 1$ and some $\rho \in (0,1]$ such that for any $(x,y) \in E_n^2$, $A \in \mathcal{E}_{n+m}$ and $n \geq 0$,

$$M_{n,n+m}(x,A) \geq \rho M_{n,n+m}(y,A),$$

where $M_{n,n+m} = (M_{n+1}\ldots M_{n+m})$ stands for the composition of the Markov kernels $M_p$ from $p = (n+1)$ to $p = (n+m)$. In this situation, following the arguments given in [5] one proves that

$$r_{n,n+m} \leq r^m/\rho \quad \text{and} \quad \beta(P_{n,n+m}) \leq (1 - r^{m-1}\rho^2)^{[(n-p)/m]}.$$

Furthermore, the constants $b(n)$ in Lemma 3.1 can be chosen such that $\bigvee_n b(n) \leq 2mr^{2m-1}/\rho^3$.

The proof of Lemma 3.1 being rather technical, it is housed in the Appendix (see Lemma A.3). One consequence of Lemma 3.1 is the following central estimate.

LEMMA 3.3. *Suppose the McKean interpretation model satisfies condition* (H) *for some finite constant* $\Gamma < \infty$. *In this situation, the martingale* $\mathcal{M}_n^N(f)$ *defined by* (8) *satisfies conditions* (H1)–(H3) *for some universal constants*

$$(a_1(n), a_2(n)) = (a_1, a_2)$$

*with the nonnegative increasing process* $C_n(f)$ *defined at* (10), *as soon as* $n \to C_n(f)$ *is strictly increasing. In addition, the constant* $a_3(n)$ *in* (H3) *can be chosen such that, for any* $n \geq 0$,

$$0 < a_3(n) \leq 4\sqrt{2}(1+\Gamma) \sup_{q=n,n-1} \sum_{p=0}^{q} r_{p,q} \beta(P_{p,q}).$$

*Furthermore, when the regularity conditions stated in Remark* 3.2 *are met for some triplet* $(m, r, \rho)$, *the constant* $a_3(n)$ *can be chosen such that* $0 < \bigvee_n a_3(n) \leq 8\sqrt{2} m r^{2m-1}(1+\Gamma)/\rho^3$.

The second step will be to get rid of the predictable term defined by (12), with the help of the following lemma:

LEMMA 3.4. *Let* $F_Z$ *be the distribution function associated to a real-valued random variable* $Z$, *and let* $W$ *be a centered Gaussian random variable with unit variance. For any pair of random variables* $(X, Y)$ *we have*

$$(20) \qquad \|F_{X+Y} - F_W\| \leq \|F_X - F_W\| + 4\mathbb{E}(|XY|) + 4\mathbb{E}(|Y|).$$

Lemma 3.4 can be proved using the Stein approach to fluctuations and it can be found, for instance, as Lemma 1.3, Chapter 11 in [8]. Since the proof of Theorem 1.3 is now a simple consequence of Lemmas 3.3 and 3.4, we postpone the proof of Lemma 3.3 and start with the

PROOF OF THEOREM 1.3. Throughout this proof, $\hat{b}(n)$ will stand for a positive constant that can change from line to line. We first notice that

$$W_{n,n}^N = \sqrt{N}(\eta_n^N - \eta_n),$$

and by (12) we have $W_{p,n}^N(f_n) = \mathcal{B}_{p,n}^N(f_n) + \mathcal{L}_{p,n}^N(f_n)$. Let us show now that the main term in the fluctuations of the c.d.f. of $W_{p,n}^N(f_n)$ is due to $\mathcal{L}_{p,n}^N(f_n)$. Indeed, it is easily checked that

$$(21) \qquad \mathbb{E}^{1/2}(|\mathcal{B}_{n,n}^N(f_n)|^2) \vee \mathbb{E}(|\mathcal{B}_{n,n}^N(f_n)|) \leq \frac{\hat{b}(n)(C_n(f))^{1/2}}{\sqrt{N}}.$$

By definition of the martingale term, it is also easily seen that

(22) $$\mathbb{E}^{1/2}(|\mathcal{L}_{n,n}^N(f_n)|^2) \leq \hat{b}(n)(C_n(f))^{1/2}.$$

Set now

$$X = \mathcal{L}_{n,n}^N(f_n)/(C_n(f))^{1/2} \quad \text{and} \quad Y = \mathcal{B}_{n,n}^N(f_n)/(C_n(f))^{1/2}.$$

The estimates (21) and (22) yield

$$\mathbb{E}(|XY|) \leq \frac{1}{C_n(f)} \mathbb{E}^{1/2}[|\mathcal{L}_{n,n}^N(f_n)|^2]\mathbb{E}^{1/2}[|\mathcal{B}_{n,n}^N(f_n)|^2] \leq \frac{\hat{b}(n)}{N^{1/2}},$$

and $\mathbb{E}(|Y|) \leq \hat{b}(n)/N^{1/2}$. Hence, applying Lemma 3.4, the proof of Theorem 1.3 can be reduced to show that

$$\sup_{u \in \mathbb{R}} \left| \mathbb{P}(\mathcal{L}_n^N(f_{n,n}) \leq u(C_n(f))^{1/2}) - \frac{1}{(2\pi)^{1/2}} \int_{-\infty}^u e^{-v^2/2}\,dv \right| \leq \frac{\hat{b}(n)}{N^{1/2}}.$$

This last estimate is now a direct consequence of Lemma 3.3 and Theorem 1.1. □

We now come to:

PROOF OF LEMMA 3.3. Let us first check that the regularity condition (H3) is satisfied. Since we have

$$\Delta \mathcal{C}_n(f) = \eta_{n-1}[K_{n,\eta_{n-1}}((f_n - K_{n,\eta_{n-1}}f_n)^2)]$$
$$= \Phi_n(\eta_{n-1})(f_n^2) - \eta_{n-1}((K_{n,\eta_{n-1}}f_n)^2),$$

we easily prove that

$$|\Delta \mathcal{C}_n^N(f) - \Delta \mathcal{C}_n(f)| \leq 2(|\Phi_n(\eta_{n-1}^N)(h_n)| + |\eta_{n-1}^N(h'_{n-1})|$$
$$+ \|K_{n,\eta_{n-1}^N}(f_n) - K_{n,\eta_{n-1}}(f_n)\|),$$

with the pair of functions $(h_n, h'_{n-1}) \in (\mathrm{Osc}(E_n) \times \mathrm{Osc}(E_{n-1}))$ defined by

$$h_n = (f_n^2 - \eta_n(f_n^2))/2,$$
$$h'_{n-1} = ((K_{n,\eta_{n-1}}f_n)^2 - \eta_{n-1}((K_{n,\eta_{n-1}}f_n)^2))/2.$$

On the other hand, under condition (H), we have that

$$\|K_{n,\eta_{n-1}^N}(f_n) - K_{n,\eta_{n-1}}(f_n)\|$$
$$\leq \int |\eta_{n-1}^N(h)|\Gamma_{n-1,f_n}(dh) + \int |\Phi_n(\eta_{n-1}^N)(h)|\Gamma'_{n-1,f_n}(dh),$$

from which we find that

$$|\Delta\mathcal{C}_n^N(f) - \Delta\mathcal{C}_n(f)|$$
$$\leq 2\left(\int |\eta_{n-1}^N(h)|\widetilde{\Gamma}_{n-1}(dh) + \mathbb{E}\left(\int |\eta_n^N(h)|\widetilde{\Gamma}'_{n-1}(dh)|\mathcal{F}_{n-1}^N\right)\right),$$

with

$$\widetilde{\Gamma}_{n-1} = \Gamma_{n-1,f_n} + \delta_{h'_{n-1}} \quad \text{and} \quad \widetilde{\Gamma}'_{n-1} = \Gamma'_{n-1,f_n} + \delta_{h_n}.$$

Applying Jensen's inequality, we get that for any $\varepsilon > 0$,

$$\mathbb{E}[e^{\varepsilon N^{1/2}|\Delta\mathcal{C}_n^N(f) - \Delta\mathcal{C}_n(f)|}]$$
$$\leq \mathbb{E}[e^{2\varepsilon N^{1/2}\{\int |\eta_{n-1}^N(h)|\widetilde{\Gamma}_{n-1}(dh) + \int |\eta_n^N(h)|\widetilde{\Gamma}'_{n-1}(dh)\}}].$$

Now, applying the Cauchy–Schwarz inequality, we obtain

$$\mathbb{E}[e^{\varepsilon N^{1/2}|\Delta\mathcal{C}_n^N(f) - \Delta\mathcal{C}_n(f)|}]$$

$$\leq \mathbb{E}^{1/2}[e^{4\varepsilon N^{1/2}\int |\eta_{n-1}^N(h)|\widetilde{\Gamma}_{n-1}(dh)}]\mathbb{E}^{1/2}[e^{4\varepsilon N^{1/2}\int |\eta_n^N(h)|\widetilde{\Gamma}'_{n-1}(dh)}].$$

If we set $\widetilde{\Gamma} = \Gamma + 1$, then using again Jensen's inequality, we find

$$\mathbb{E}[e^{\varepsilon N^{1/2}|\Delta\mathcal{C}_n^N(f) - \Delta\mathcal{C}_n(f)|}]^2 \leq \int \mathbb{E}[e^{4\varepsilon N^{1/2}\widetilde{\Gamma}|\eta_{n-1}^N(h)|}]\frac{\widetilde{\Gamma}_{n-1}(dh)}{\widetilde{\Gamma}_{n-1}(1)}$$
$$\times \int \mathbb{E}[e^{4\varepsilon N^{1/2}\widetilde{\Gamma}|\eta_n^N(h)|}]\frac{\widetilde{\Gamma}'_{n-1}(dh)}{\widetilde{\Gamma}'_{n-1}(1)},$$

from which we get

$$\mathbb{E}[e^{\varepsilon N^{1/2}|\Delta\mathcal{C}_n^N(f) - \Delta\mathcal{C}_n(f)|}] \leq \sup_{h \in \text{Osc}(E_p), p=n,n-1} \mathbb{E}(e^{4\varepsilon N^{1/2}\widetilde{\Gamma}|\eta_p^N(h)|}).$$

Using Lemma A.3, we conclude that

$$\mathbb{E}[e^{\varepsilon N^{1/2}|\Delta\mathcal{C}_n^N(f) - \Delta\mathcal{C}_n(f)|}] \leq (1 + \varepsilon a_3(n))e^{\varepsilon^2 a_3^2(n)}$$

for some finite constant $a_3(n)$ such that

$$a_3(n) \leq 4\sqrt{2}\widetilde{\Gamma} \sup_{q=n,n-1} \sum_{p=0}^{q} r_{p,q}\beta(P_{p,q}).$$

To prove that (H2) is met, we first recall that

(23) $$|\mathbb{E}[e^{i\lambda N^{1/2}\mathcal{M}_n^N(f)}]| \leq \mathbb{E}[|\mathbb{E}[e^{i\lambda N^{1/2}\Delta\mathcal{M}_n^N(f)}|F_{n-1}^N]|].$$

Then we use a standard symmetrization technique: given the particle model $\xi_p$ up to time $p \leq n-1$, we let $\bar{\eta}_n^N$ be an auxiliary independent copy of $\eta_n^N$. In other words, $\bar{\eta}_n^N$ is the empirical measure associated to an independent copy $\bar{\xi}_n$ of the configuration of the system $\xi_n$ at time $n$. With some obvious abuse of notation, we readily check that

$$|\mathbb{E}[e^{i\lambda N^{1/2}\Delta\mathcal{M}_n^N(f)}|F_{n-1}^N]|^2 = \mathbb{E}[e^{i\lambda N^{1/2}[\Delta\mathcal{M}_n^N(f)-\Delta\overline{\mathcal{M}}_n^N(f)]}|F_{n-1}^N],$$

where $\Delta\overline{\mathcal{M}}_n^N(f) = [\bar{\eta}_n^N(f_n) - \Phi_n(\eta_{n-1}^N)(f_n)]$. We deduce from this that

$$|\mathbb{E}[e^{i\lambda N^{1/2}\Delta\mathcal{M}_n^N(f)}|F_{n-1}^N]|^2 = \prod_{j=1}^N \mathbb{E}[e^{i(\lambda/N^{1/2})[f_n(\xi_n^j)-f_n(\bar{\xi}_n^j)]}|\mathcal{F}_{n-1}^N].$$

Since the random variables $[f_n(\xi_n^j) - f_n(\bar{\xi}_n^j)]$ and $-[f_n(\xi_n^j) - f_n(\bar{\xi}_n^j)]$ have the same law, their characteristic functions are real, and we have

$$\mathbb{E}[e^{i(\lambda/N^{1/2})[f_n(\xi_n^j)-f_n(\bar{\xi}_n^j)]}|\mathcal{F}_{n-1}^N] = \mathbb{E}\left[\cos\left(\frac{\lambda}{N^{1/2}}[f_n(\xi_n^j) - f_n(\bar{\xi}_n^j)]\right)\bigg|\mathcal{F}_{n-1}^N\right].$$

Using now the elementary inequalities

$$\cos u \leq 1 - u^2/2 + |u|^3/3!, \qquad 1 + u \leq e^u, \qquad |u-v|^3 \leq 4(|u|^3 + |v|^3),$$

we get that

$$\mathbb{E}[e^{i(\lambda/N^{1/2})[f_n(\xi_n^j)-f_n(\bar{\xi}_n^j)]}|\mathcal{F}_{n-1}^N]$$

$$\leq 1 - \frac{\lambda^2}{N}K_{n,\eta_{n-1}^N}(f_n - K_{n,\eta_{n-1}^N}(f_n))^2(\xi_n^j) + \frac{c\lambda^3}{N^{3/2}}$$

$$\leq e^{-(\lambda^2/N)}K_{n,\eta_{n-1}^N}(f_n - K_{n,\eta_{n-1}^N}(f_n))^2(\xi_n^j) + \frac{c\lambda^3}{N^{3/2}}.$$

Multiplying over $j$, we obtain

$$|\mathbb{E}[e^{i\lambda N^{1/2}\Delta\mathcal{M}_n^N(f)}|\mathcal{F}_{n-1}^N]|^2 \leq e^{-\lambda^2 \Delta\mathcal{C}_n^N(f) + c\lambda^3/N^{1/2}},$$

and by (23) we conclude that condition (H2) is met with $a_2(n) = c/2$.

We now come to the proof of (H1). By definition of the particle model associated to a given collection of transitions $K_{n,\eta}$ we have

$$\mathbb{E}[e^{i\lambda N^{1/2}\Delta\mathcal{M}_n^N(f)+(\lambda^2/2)\Delta\mathcal{C}_n^N(f)}|\mathcal{F}_{n-1}^N]$$

$$= \prod_{j=1}^N [K_{n,\eta_{n-1}^N}(e^{i(\lambda/N^{1/2})\tilde{f}_n^j+(\lambda^2/(2N))\Delta\mathcal{C}_n^N(f)})](\xi_{n-1}^j),$$

with the random function $\tilde{f}_n^j = (f_n - K_{n,\eta_{n-1}^N}(f_n))(\xi_{n-1}^j)$. Using the elementary inequality

$$\left|e^z - \left(1 + z + \frac{z^2}{2}\right)\right| \leq e^{|z|}\frac{|z|^3}{3!},$$

it is easily seen that, for any $\lambda \leq N^{1/2}$, we have

$$e^{i(\lambda/N^{1/2})\tilde{f}_n^j + (\lambda^2/(2N))\Delta\mathcal{C}_n^N(f)}$$

$$= 1 + i\frac{\lambda}{N^{1/2}}\tilde{f}_n^j + \frac{\lambda^2}{2N}[\Delta\mathcal{C}_n^N(f) - (\tilde{f}_n^j)^2] + r_{n,1}^N(f),$$

with $|r_{n,1}^N(f)| \leq c\lambda^3 N^{-3/2}$. This clearly implies that, for any $\lambda \leq N^{1/2}$,

$$[K_{n,\eta_{n-1}^N}(e^{i(\lambda/N^{1/2})\tilde{f}_n^j + (\lambda^2/(2N))\Delta\mathcal{C}_n^N(f)})](\xi_{n-1}^j)$$

$$= 1 + \frac{\lambda^2}{2N}[\Delta\mathcal{C}_n^N(f) - K_{n,\eta_{n-1}^N}(\tilde{f}_n^j)^2(\xi_{n-1}^j)] + r_{n,2}^N(f),$$

with $|r_{n,2}^N(f)| \leq c\lambda^3 N^{-3/2}$. It is now convenient to notice that for any $\lambda \leq N^{1/2}$,

$$\left|\frac{\lambda^2}{2N}[\Delta\mathcal{C}_n^N(f) - K_{n,\eta_{n-1}^N}(\tilde{f}_n^j)^2(\xi_{n-1}^j)] + r_{n,2}^N(f)\right| \leq \frac{c\lambda}{N^{1/2}}.$$

On the other hand, for any $|z| \leq 1/2$ and with the principal value of the logarithm we recall that

$$\log(1+z) = z - \int_0^z \frac{u}{1+u}\,du = z - z^2\int_0^1 \frac{t}{1+tz}\,dt.$$

Since for any $|z| \leq 1/2$ and $t \in [0,1]$ we have $|1+tz| \geq 1/2$, we find that for any $|z| \leq 1/2$ we have $|\log(1+z) - z| \leq |z^2|$. The previous computations show that there exists some universal constant $c_0 \in (0,1)$ such that for any $\lambda \leq c_0 N^{1/2}$ we have

$$\log K_{n,\eta_{n-1}^N}(e^{i(\lambda/N^{1/2})\tilde{f}_n^j + (\lambda^2/(2N))\Delta\mathcal{C}_n^N(f)})(\xi_{n-1}^j)$$

$$= \frac{\lambda^2}{2N}[\Delta\mathcal{C}_n^N(f) - K_{n,\eta_{n-1}^N}(\tilde{f}_n^j)^2(\xi_{n-1}^j)] + r_{n,3}^N(f),$$

with $|r_{n,3}^N(f)| \leq c\lambda^3/N^{3/2}$. Summing over $j$, we see that for any $\lambda \leq c_0 N^{1/2}$,

$$\left|\sum_{j=1}^N \log K_{n,\eta_{n-1}^N}(e^{i(\lambda/N^{1/2})\tilde{f}_n^j + (\lambda^2/(2N))\Delta\mathcal{C}_n^N(f)})(\xi_{n-1}^j)\right| \leq c\lambda^3/N^{1/2}.$$

Finally, using the elementary inequality $|e^z - 1| \leq |z|e^{|z|}$, we conclude that, for any $\lambda \leq c_0 N^{1/2}$,

$$|\mathbb{E}[e^{i\lambda N^{1/2}\Delta\mathcal{M}_n^N(f) + (\lambda^2/2)\Delta\mathcal{C}_n^N(f)}|\mathcal{F}_{n-1}^N] - 1| \leq c\frac{\lambda^3}{N^{1/2}}e^{c\lambda^3/N^{1/2}}.$$

This readily implies that there exists some universal positive constant $c_1$ such that, for any $\lambda^3 \leq c_1 N^{1/2}$, we have

$$|\mathbb{E}[e^{i\lambda N^{1/2}\Delta\mathcal{M}_n^N(f) + (\lambda^2/2)\Delta\mathcal{C}_n^N(f)}|\mathcal{F}_{n-1}^N] - 1| \leq c\frac{\lambda^3}{N^{1/2}},$$

which proves that condition (H1) is met with $a_1(n) = c$ and $c_1(n) = 1$. □

## APPENDIX

### A.1. Doob type decompositions.

PROPOSITION A.1 ([5]). *Let $(\overline{Q}_{p,n})_{p \leq n}$ be the semigroup defined at* (11). *For $f_n \in \mathcal{B}_b(E_n)$ and $p \leq n$ we set $f_{p,n} = \overline{Q}_{p,n}(f_n - \eta_n f_n)$. Then we have the following decomposition:*

$$\eta_p^N(f_{p,n}) = \mathcal{A}_{p,n}^N(f_n) + \mathcal{M}_{p,n}^N(f_n),$$

*with the predictable and martingale sequences $\mathcal{A}_{p,n}^N(f_n)$ and $\mathcal{M}_{p,n}^N(f_n)$ given by*

$$(\text{A.1}) \qquad \mathcal{A}_{p,n}^N(f_n) = \sum_{q=1}^{p} [1 - \eta_{q-1}^N(\overline{Q}_{q-1,q}1)]\Phi_q(\eta_{q-1}^N)(f_{q,n}),$$

$$(\text{A.2}) \qquad \mathcal{M}_{p,n}^N(f_n) = \sum_{q=0}^{p} [\eta_q^N(f_{q,n}) - \Phi_q(\eta_{q-1}^N)f_{q,n}],$$

*with the usual convention $\Phi_0(\eta_{-1}^N) = \eta_0$.*

PROOF. Note that for any $\varphi_n \in \mathcal{B}_b(E_n)$ we have the decomposition

$$\eta_p^N(\overline{Q}_{p,n}\varphi_n) - \eta_0^N(\overline{Q}_{0,n}\varphi_n) = \sum_{q=1}^{p} \bar{\delta}_q,$$

with $\bar{\delta}_q = \eta_q^N(\overline{Q}_{q,n}\varphi_n) - \eta_{q-1}^N(\overline{Q}_{q-1,n}\varphi_n)$. Choose now $\varphi_n = f_n - \eta_n f_n$. For $q \leq p$, we have, by definition of $f_{q,n}$,

$$\bar{\delta}_q = \eta_q^N(f_{q,n}) - \eta_{q-1}^N(\overline{Q}_{q-1,n}\varphi_n) = U_1 + U_2,$$

with

$$U_1 = \eta_q^N(f_{q,n}) - [\Phi_q(\eta_{q-1}^N)](f_{q,n}),$$
$$U_2 = [\Phi_q(\eta_{q-1}^N)](f_{q,n}) - \eta_{q-1}^N(\overline{Q}_{q-1,n}\varphi_n).$$

In order to show (12), it is thus enough to verify that

$$(\text{A.3}) \qquad \eta_{q-1}^N(\overline{Q}_{q-1,n}\varphi_n) = \eta_{q-1}^N(\overline{Q}_{q-1,q}1)[\Phi_q(\eta_{q-1}^N)](f_{q,n}).$$

However, we have

$$\overline{Q}_{q-1,n}\varphi_n = \frac{\gamma_{q-1}(1)}{\gamma_n(1)}Q_{q-1,n}\varphi_n = \frac{\gamma_{q-1}(1)}{\gamma_n(1)}Q_q(Q_{q,n}\varphi_n)$$
$$= \frac{\gamma_{q-1}(1)}{\gamma_q(1)}Q_q(\overline{Q}_{q,n}\varphi_n) = \frac{\gamma_{q-1}(1)}{\gamma_q(1)}Q_q(f_{q,n}),$$

and hence

$$\eta_{q-1}^N(\overline{Q}_{q-1,n}\varphi_n) = \frac{\gamma_{q-1}(1)}{\gamma_q(1)}\eta_{q-1}^N(Q_q(f_{q,n}))$$

(A.4)
$$= \frac{\gamma_{q-1}(1)\eta_{q-1}^N(G_{q-1})}{\gamma_q(1)}[\Phi_q(\eta_{q-1}^N)](f_{q,n}).$$

On the other hand, for $q \geq 1$ and $x_{q-1} \in E_{q-1}$,

$$[\overline{Q}_{q-1,q}(1)](x_{q-1}) = \frac{\gamma_{q-1}(1)}{\gamma_q(1)}[Q_{q-1,q}(1)](x_{q-1}),$$

and

$$[Q_{q-1,q}(1)](x_{q-1}) = \int_{E_q} G_{q-1}(x_{q-1})M_q(x_{q-1},dx_q) = G_{q-1}(x_{q-1}),$$

which yields

$$\frac{\gamma_{q-1}(1)G_{q-1}}{\gamma_q(1)} = \overline{Q}_{q-1,q}(1).$$

Plugging this last equality into (A.4), we get (A.3), and hence (12). The martingale property of $\mathcal{M}_{p,n}^N(f_n)$ is readily checked. □

**A.2. Some asymptotic estimates.** The next lemma provides a refined version of Burkholder type inequalities for independent sequences of random variables.

LEMMA A.2 ([6]). *Let $m(X) = \frac{1}{N}\sum_{i=1}^N \delta_{X^i}$ be the $N$-empirical measure associated to a collection of independent random variables $X^i$, with respective distributions $\mu_i$ on some measurable space $(E,\mathcal{E})$. For any sequence of $\mathcal{E}$-measurable functions $h_i$ such that $\mu_i(h_i) = 0$ and $\sigma^2(h) = \frac{1}{N}\sum_{i=1}^N \operatorname{osc}^2(h_i) < \infty$, we have for any integer $p \geq 1$*

(A.5) $$\sqrt{N}\mathbb{E}(|m(X)(h)|^p)^{1/p} \leq d(p)^{1/p}\sigma(h),$$

*with the sequence of finite constants $(d(n))_{n\geq 0}$ defined for any $n \geq 1$ by the formulae*

(A.6) $$d(2n) = (2n)_n 2^{-n} \quad \text{and} \quad d(2n-1) = \frac{(2n-1)_n}{\sqrt{n-1/2}}2^{-(n-1/2)}.$$

The extension of Lemma A.2 to the interacting particle measures $\eta_n^N$ and the Feynman–Kac flow $\eta_n$ is the following

LEMMA A.3. *Let $(d(p))_{p \geq 1}$ be the sequence introduced in* (A.6). *For any McKean interpretation model and for any $n \geq 0$, $p \geq 1$, $f_n \in \mathrm{Osc}(E_n)$ and $\varepsilon > 0$, we have*

$$\mathbb{E}(|[\eta_n^N - \eta_n](f_n)|^p)^{1/p} \leq d(p)^{1/p} b(n)/\sqrt{N},$$

$$\mathbb{E}[e^{\varepsilon N^{1/2} |\eta_n^N(f_n) - \eta_n(f_n)|}] \leq (1 + \varepsilon 2^{-1/2} b(n)) e^{(\varepsilon b(n))^2/2},$$

*for some finite constant $b(n)$ such that $b(n) \leq 2 \sum_{q=0}^n r_{q,n} \beta(P_{q,n})$.*

PROOF. The proof is based on the following decomposition:

$$(A.7) \qquad \eta_n^N - \eta_n = \sum_{q=0}^n [\Phi_{q,n}(\eta_q^N) - \Phi_{q,n}(\Phi_q(\eta_{q-1}^N))].$$

We introduce the random potential functions

$$G_{q,n}^N : x_q \in E_q \to G_{q,n}^N(x_q) = \frac{G_{q,n}}{\Phi_q(\eta_{q-1}^N)(G_{q,n})} \in (0, \infty)$$

and the random bounded operators $P_{q,n}^N$ from $\mathcal{B}_b(E_n)$ into $\mathcal{B}_b(E_q)$ defined for any $(f_n, x_q) \in (\mathcal{B}_b(E_n) \times E_q)$ by

$$P_{q,n}^N(f_n)(x_q) = P_{q,n}(f_n - \Phi_{q,n}(\Phi_q(\eta_{q-1}^N))(f_n))(x_q)$$

$$= \int (P_{q,n} f(x_q) - P_{q,n} f(y_q)) G_{q,n}^N(y_q) \Phi_q(\eta_{q-1}^N)(dy_q).$$

We associate to the pair $(G_{q,n}^N, P_{q,n}^N)$ the random bounded and integral operator $Q_{q,n}^N$ from $\mathcal{B}_b(E_n)$ into $\mathcal{B}_b(E_q)$ defined for any $(f_n, x_q) \in (\mathcal{B}_b(E_n) \times E_q)$ by

$$Q_{q,n}^N(f_n)(x_q) = G_{q,n}^N(x_q) \times P_{q,n}^N(f_n)(x_q).$$

Each "local" term in (A.7) can be expressed in terms of $Q_{q,n}^N$ as follows. For any $q \leq n$ and $f_n \in \mathcal{B}_b(E_n)$ with $\mathrm{osc}(f_n) \leq 1$ we have

$$\Phi_{q,n}(\eta_q^N)([f_n - \Phi_{q,n}(\Phi_q(\eta_{q-1}^N))(f_n)])$$

$$= \frac{1}{\eta_q^N(G_{q,n})} \eta_q^N(G_{q,n} P_{q,n}[f_n - \Phi_{q,n}(\Phi_q(\eta_{q-1}^N))(f_n)])$$

$$= \frac{\eta_q^N Q_{q,n}^N(f_n)}{\eta_q^N Q_{q,n}^N(1)}.$$

By construction we also observe that

$$\Phi_q(\eta_{q-1}^N)(G_{q,n}^N) = 1 \quad \text{and} \quad \Phi_q(\eta_{q-1}^N)(Q_{q,n}^N(f_n)) = 0.$$

The above considerations easily yield the decomposition

$$\Phi_{q,n}(\eta_q^N) - \Phi_{q,n}(\Phi_q(\eta_{q-1}^N)) = \frac{1}{\eta_q^N(G_{q,n}^N)}[\eta_q^N - \Phi_q(\eta_{q-1}^N)]Q_{q,n}^N.$$

Using the properties of the Dobrushin contraction coefficient, we also have

$$\|P_{q,n}^N(f_n)\| \leq \mathrm{osc}(P_{q,n}f) \leq \beta(P_{q,n}),$$
$$\|Q_{q,n}^N(f_n)\| \leq \|G_{q,n}^N\|\|P_{q,n}^N(f_n)\| \leq \|G_{q,n}^N\|\beta(P_{q,n}),$$

and from these estimates, we readily prove the inequality

$$|[\Phi_{q,n}(\eta_q^N) - \Phi_{q,n}(\Phi_q(\eta_{q-1}^N))](f_n)|$$
$$\leq r_{q,n}\beta(P_{q,n})|[\eta_q^N - \Phi_q(\eta_{q-1}^N)]\overline{Q}_{q,n}^N(f_n)|,$$

with $\overline{Q}_{q,n}^N(f_n) = Q_{q,n}^N(f_n)/\|Q_{q,n}^N(f_n)\|$. Now, using Lemma A.2, we check that for any $p \geq 1$ we have

$$\sqrt{N}\mathbb{E}(|[\eta_q^N - \Phi_q(\eta_{q-1}^N)]\overline{Q}_{q,n}^N(f_n)|^p|F_{q-1}^N)^{1/p} \leq 2d(p)^{1/p},$$

with the sequence of finite constants $d(p)$ introduced in (A.6). This ends the proof of the first assertion. The $\mathbb{L}_n$-inequalities stated in Lemma A.2 clearly implies that, for any $\varepsilon > 0$,

$$\mathbb{E}(e^{\varepsilon|\eta_n^N(f_n) - \eta_n(f_n)|}) = \sum_{n \geq 0} \frac{\varepsilon^{2n}}{(2n)!}\mathbb{E}(|\eta_n^N(f_n) - \eta_n(f_n)|^{2n})$$
$$+ \sum_{n \geq 0} \frac{\varepsilon^{2n+1}}{(2n+1)!}\mathbb{E}(|\eta_n^N(f_n) - \eta_n(f_n)|^{2n+1})$$
$$\leq \sum_{n \geq 0} \frac{1}{n!}\left(\frac{\varepsilon^2 b(n)^2}{2N}\right)^n + \sum_{n \geq 0} \frac{1}{n!}\left(\frac{\varepsilon^2 b(n)^2}{2N}\right)^{n+1/2},$$

from which we conclude that

$$\mathbb{E}(e^{\varepsilon|\eta_n^N(f_n) - \eta_n(f_n)|}) = \left(1 + \frac{\varepsilon b(n)}{\sqrt{2N}}\right) \sum_{n \geq 0} \frac{1}{n!}\left(\frac{\varepsilon^2 b(n)^2}{2N}\right)^n$$
$$= \left(1 + \frac{\varepsilon b(n)}{\sqrt{2N}}\right) e^{(\varepsilon^2/2N)b(n)^2}.$$

We end the proof of the lemma by replacing $\varepsilon$ by $\varepsilon\sqrt{N}$. $\square$

## REFERENCES

[1] DEL MORAL, P. (1998). Measure valued processes and interacting particle systems. Application to nonlinear filtering problems. *Ann. Appl. Probab.* **8** 438–495. MR1624949
[2] DEL MORAL, P. and GUIONNET, A. (1999). Central limit theorem for nonlinear filtering and interacting particle systems. *Ann. Appl. Probab.* **9** 275–297. MR1687359
[3] DEL MORAL, P. and JACOD, J. (2002). The Monte Carlo method for filtering with discrete time observations. Central limit theorems. In *Numerical Methods and Stochastics* (T. J. Lyons and T. S. Salisbury, eds.) 29–53. Amer. Math. Soc., Providence, RI.
[4] DEL MORAL, P. and LEDOUX, M. (2000). Convergence of empirical processes for interacting particle systems with applications to nonlinear filtering. *J. Theoret. Probab.* **13** 225–257. MR1744985
[5] DEL MORAL, P. and MICLO, L. (2000). Branching and interacting particle systems approximations of Feynman–Kac formulae with applications to non-linear filtering. *Séminaire de Probabilités XXXIV. Lecture Notes in Math.* **1729** 1–145. Springer, New York. MR1768060
[6] DEL MORAL, P., MICLO, L. and VIENS, F. (2003). Precise propagations of chaos estimates for Feynman–Kac and genealogical particle models. Technical Report 03-01, Center for Statistical Decision Sciences and Dept. Statistics, Purdue Univ.
[7] DYNKIN, E. B. and MANDELBAUM, A. (1983). Symmetric statistics, Poisson processes and multiple Wiener integrals. *Ann. Statist.* **11** 739–745. MR707925
[8] SHORACK, G. R. (2000). *Probability for Statisticians.* Springer, New York. MR1762415
[9] SHIGA, T. and TANAKA, H. (1985). Central limit theorem for a system of Markovian particles with mean field interaction. *Z. Wahrsch. Verw. Gebiete* **69** 439–459. MR787607

LABORATOIRE J.-A. DIEUDONNEÉ
U.M.R. 6621 DU C.N.R.S.
UNIVERSITÉ DE NICE-SOPHIA ANTIPOLIS
PARC VALROSE
06108 NICE CEDEX 02
FRANCE
E-MAIL: delmoral@math.unice.fr

INSTITUT ÉLIE CARTAN
UNIVERSITÉ HENRI POINCARÉ NANCY 1
BP 239, 54506-VANDOEUVRE-LÈS-NANCY
FRANCE
E-MAIL: tindel@iecn.u-nancy.fr